\documentclass{article}
\usepackage[utf8]{inputenc}
\usepackage{amssymb,amsmath,latexsym}
\usepackage[english]{babel}
\usepackage{amscd,amsxtra}
\usepackage{amsthm}
\usepackage{verbatim}
\usepackage{graphicx}
\usepackage{enumitem}
\usepackage{color}
\usepackage{tikz}
\usepackage{caption}
\usepackage[rightcaption]{sidecap}
\usepackage[colorinlistoftodos]{todonotes}
\usetikzlibrary{arrows.meta}
\tikzstyle{point}=[ball color=white, circle, draw=black, inner sep=0.1cm]
\usepackage{authblk}

\newtheorem{theorem}{Theorem}[section]

\newtheorem{lemma}[theorem]{Lemma}

\newtheorem{proposition}[theorem]{Proposition}
\newtheorem{remark}[theorem]{Remark}

\title{Graph multicoloring in solving node malfunction and attacks in a secret-sharing network}
\author[1]{Tanja Vojković \thanks{Corresponding author: tanja@pmfst.hr, Rudjera Bo\v skovi\' ca 33, 21000 Split, Croatia}}
\author[1]{Damir Vukičević}
\affil[1]{Faculty of Science, University of Split}
\date{May 2023}

\begin{document}
\maketitle
\begin{abstract}
   We observe a network scenario where parts of a secret are distributed among its nodes. Within the network, a group of attackers is actively trying to obtain the complete secret, while there is also the issue of some nodes malfunctioning or being absent. In this paper, we address this problem by employing graph multicoloring techniques, focusing on the case of a single attacker and varying numbers of malfunctioning nodes.
\end{abstract}

Keywords: graph multicoloring, secret-sharing, networks

\section{Introduction}
This paper builds upon our previous work titled 'Multicoloring of Graphs to Secure a Secret' \cite{nas}, where we introduced a specific type of attack on secret-sharing networks and proposed a solution to protect secrets in such scenarios. In a secret-sharing network, a group divides a secret into parts and distributes them among participants, represented as nodes in a network of common acquaintances. Attackers infiltrate this network with dual objectives: obtaining all parts of the secret and disrupting network connections to prevent secret reconstruction. Our research seeks to answer questions such as how many secret pieces are required to ensure its safety against a specified number of attackers, and what network topology and piece distribution configurations achieve this security. 
We model the secret pieces with colors and observe the problem as the problem of graph multicoloring. In our previous work \cite{nas}, we tackled the problem of securing secrets for 1, 2, 3, and 4 attackers. In this paper, we extend our analysis to address the issue of network nodes malfunctioning.\\
In real-world networks, node failures are common, and designing network structures must account for this. If the participants (nodes) in a network are persons guarding a piece of the secret one can imagine a node failure as a person simply giving up on the task given or being removed from the network by some outside means and if the nodes represent some piece of technology then a failure can mean a simple malfunction. In Section \ref{def} we give definitions of basic concepts and preliminaries needed to study the problem, while Section \ref{main} provides our main results.

%-----------------------------------------------------------------------------------------------

\section{Definitions and preliminary results}\label{def}
We will use standard definitions and notation from graph theory, \cite{gross}. By $N(u)=N_G(u)$ we denote the set of neighbors of vertex $u$ in the graph $G$, and for $A\subseteq V(G)$, by $N_G(A)$ we denote the set of neighbors of all the vertices in $A$. Further we denote,
$$N_0(u)=N(u)\cup \{u\},$$
$$N_0(A)=N(A)\cup A.$$

By $G-u$ we denote a graph obtained from graph $G$ by deleting the vertex $u\in V(G)$ and all its incident edges, and for $A\subseteq V(G)$, by $G-A$, we denote a graph obtained from graph $G$ by deleting all the vertices in $A$ and their incident edges.
Let us first define a graph multicoloring, and the specific multicoloring we used to solve the problem in \cite{nas}.\\

Let $k\in\mathbb{N}$. \textbf{Vertex $k$-multicoloring} of a graph $G$ is a function $\kappa:V(G)\to \mathcal{P}(\{1,...,k\})$. So, while in vertex coloring of a graph, each vertex is colored with one color, in multicoloring, each vertex is colored by some subset of the set of $k$ colors. Let $a\in\mathbb{N}$. 
A vertex $k$-multicoloring of a graph $G$ is called an \textbf{$a$-resistant vertex $k$-multicoloring} if the following holds:

For each $A\subseteq V(G)$ with $a$ vertices, there is a component $H$ of the graph $G-N_0(A)$ such that
$${\displaystyle\bigcup\limits_{u\in V(H)}}\kappa(u)=\{1,...,k\}\text{.}$$

An $a$-resistant vertex $k$-multicoloring is called a \textbf{highly $a$-resistant vertex $k$-multicoloring} if for each $A\subseteq V(G)$ with $a$ vertices it holds that

\begin{center}
	${\displaystyle\bigcup\limits_{u\in A}}\kappa(u)\neq\{1,...,k\}$.
\end{center}

We can see that the definition of an $a$-resistant vertex $k$-multicoloring refers to any $a$ attackers not being able to disrupt the secret-sharing in the graph, because after the $a$ attackers and their neighbors leave the graph, there will exist a component $H$ of the remaining graph, which has all the colors, i.e. all parts of the secret. The condition that $a$ attackers do not collect all $k$ pieces of the secret is required for a multicoloring to be a highly $a$-resistant vertex $k$-multicoloring.
We will denote by $a-HR$ the condition that 
${\displaystyle\bigcup\limits_{u\in A}}\kappa(u)\neq\{1,...,k\}$, for each $A\subseteq V(G)$ with $a$ vertices.\\ 

Now let us take into account possible missing or malfunctioning nodes in the network. Let $a,m\in\mathbb{N}_0$. A vertex $k$-multicoloring of a graph $G$ is called an \textbf{$(a,m)$-resistant vertex $k$-multicoloring} if the following holds:

For all subsets $A,M\subseteq V(G)$, such that $|A|=a$, $|M|=m$,  there exists a component $H$ of the graph $G-(N_{0}(A)\cup M)$ such that
$${\displaystyle\bigcup\limits_{u\in V(H)}}\kappa(u)=\{1,...,k\}\text{.}$$ 
An $(a,m)$-resistant vertex $k$-multicoloring is called a \textbf{highly $(a,m)$-resistant vertex $k$-multicoloring} if the $a-HR$ condition holds.\\

%_-------------------------------------- mozda ne
Let $a,m\in\mathbb{N}_0$. If there exists a graph $G$ with $n$ vertices, and a $k\in\mathbb{N}$, such that $G$ has a highly $(a,m)$-resistant vertex $k$-multicoloring $\kappa$, we say that the pair $(G,\kappa)$ \textbf{realizes} the quadruple $(a,m,n,k)$. Moreover, we will often say that $G$ \textbf{admits} a highly $(a,m)$-resistant vertex $k$-multicoloring $\kappa$.
%-----------------------------------------------------

It is easy to see that an analogous result to the Theorem 3.1. from \cite{nas} holds, which we state as a proposition and omit the proof.

\begin{proposition}\label{uvodna}
Let $(G,\kappa)$ be a pair that realizes a quadruple $(a,m,n,k)$, for some $a,m,n,k\in\mathbb N_0$.
\begin{enumerate}
    \item There exists a graph $G'$ with $n+1$ vertices and a highly $(a,m)$-resistant vertex $k$-multicoloring $\kappa'$ of $G'$ such that $(G',\kappa')$ realizes the quadruple $(a,m,n+1,k)$.
    \item There exists a highly $(a,m)$-resistant vertex $(k+1)$-multicoloring $\kappa'$ of $G$. 
\end{enumerate}
\end{proposition}

On the basis of this proposition, we can search for the minimal number of colors $k$, for a highly $(a,m)$-resistant vertex $k$-multicoloring to exist, for given $a,m\in\mathbb{N}_0$, and for a graph $G$ with $n$ vertices that admits such multicoloring.\\

Let $a,m\in\mathbb{N}_0$. By $K(a,m,n)$ we denote the minimal number of colors such that there exists a graph $G$ with $n$ vertices and a highly $(a,m)$-resistant vertex $K(a,m,n)$-multicoloring of $G$. If, for some $a,m,n\in\mathbb{N}_0$, a graph with $n$ vertices that admits a highly $(a,m)$-resistant vertex $k$-multicoloring doesn't exist, for any $k\in\mathbb{N}_0$, we write $K(a,m,n)=\infty$.\\
Note that in \cite{nas} we found $K(a,0,n)$, for $a=1,2,3,4$.

%-------------------------------------------------------------------------

\section{Main results}\label{main}

The problem of finding the number $K(a,m,n)$ can be studied for any given $a$ and $m$. In this paper we solve the case of $a=1$ and $m\in\mathbb{N}$. 

It is easy to see that for $a=0$ the following holds.

\begin{proposition}
    Let $m\in\mathbb{N}$.
    $$K(0,m,n)=\left\{
\begin{array}
[c]{cc}%
+\infty, & n\leq m;\\
1, & n>m.
\end{array}
\right.$$
\end{proposition}

For the proof of Theorem \ref{glavni} we first need two lemmas.

%---------------------------------------------------------------------------------

\begin{lemma}\label{prosirenje}
    Let $G$ be a graph with $n$ vertices, and $f,f':V(G)\rightarrow\mathcal{P}(\{1,2,...,k\})$, vertex multicoloring functions on $G$ for which $1$-HR condition holds, such that for each $u\in V(G)$: 
    \begin{enumerate}
        \item $|f(u)|=k-1$;
        \item $f'(u)\subseteq f(u)$.
    \end{enumerate}
    If $(G,f)$ does not realize a quadruple $(1,m,n,k)$, for some $m\in\mathbb{N}$, then neither does $(G,f')$.
\end{lemma}
\begin{proof}
     Function $f$ is obtained as an extension of $f'$ such that to each vertex $u\in V(G)$, for which $|f'(u)|\leq k-2$ holds, we add some colors until they have all but one color. It is easy to see that if $f$ is not a highly $(1,m)$-resistant $k$-multicoloring of $G$, then neither if $f'$.
\end{proof}

%------------------------------------------------------------------------------

\begin{lemma}\label{pad}
Let $m\in\mathbb{N}$ and let $f:\mathbb{R}	\rightarrow\mathbb{R}$ be a function defined with	
	$$f(x)=x+\left\lfloor \frac{m}{x}\right\rfloor \text{.}$$
	Let $x_{m}\in\mathbb{R}$ such that $f(x_{m})\leq f(x)$, for all $x\in\mathbb{R}$, and $x_{m}$ is the smallest value for which that holds. It holds
	$f|_{\left\langle -\infty,x_{m}\right\rangle }$ is monotonically decreasing.
\end{lemma}
\begin{proof}
Let $m=b^{2}+c$, $0\leq c\leq2b$. Let us prove that $x_{m}\in\{b,b+1,b+2\}$.
	
	1) Let $x=b$. We have
	
	$$x+\left\lfloor \dfrac{m}{x}\right\rfloor =b+\left\lfloor \dfrac{b^{2}+c}%
	{b}\right\rfloor =2b+\left\lfloor \frac{c}{b}\right\rfloor \text{.}$$

	2) For $x=b-l$, $l\in\mathbb{R}$ we have
	$$
	x+\left\lfloor \dfrac{m}{x}\right\rfloor =b-l+\left\lfloor \dfrac{b^{2}%
		+c}{b-l}\right\rfloor=$$
	$$=b-l+\left\lfloor \dfrac{b^{2}-l^{2}+l^{2}+c}%
	{b-l}\right\rfloor =2b+\left\lfloor \frac{l^{2}+c}{b-l}\right\rfloor \text{.}%
	$$
	This expression is obviously larger than the expression from 1).
	
	3) For $x=b+l$, $l\in\mathbb{R}$ we have
	$$
	x+\left\lfloor \dfrac{m}{x}\right\rfloor=b+l+\left\lfloor \dfrac{b^{2}%
		+c}{b+l}\right\rfloor =b+l+\left\lfloor \dfrac{b^{2}-l^{2}+l^{2}+c}%
	{b+l}\right\rfloor =$$
	$$=2b+\left\lfloor \frac{l^{2}+c}{b+l}\right\rfloor \text{.}%
	$$
	By further calculation we can easily see that this expression is smaller or
	equal to the expression from 1) for $0\leq l\leq2$.
	
	This proves $x_{m}\in\{b,b+1,b+2\}$. Now let $y_{1},y_{2}\in\mathbb{R}$ such that $y_{1}\leq y_{2}\leq x_{m}$. Let us prove $f(y_{1})\geq f(y_{2})$.
	We consider some subcases depending on $x_{m}$.
	
	a) $x_{m}=b$. Let $y_{1}+l_{1}=b$, $y_{2}+l_{2}=b$, $l_{1},l_{2}\in\mathbb{R}$. It holds $l_{1}\geq l_{2}$. We have
	$$
	f(y_{1})=y_{1}+\left\lfloor \frac{m}{y_{1}}\right\rfloor =b-l_{1}%
	+\left\lfloor \frac{m}{b-l_{1}}\right\rfloor =$$
	$$=b-l_{1}+\left\lfloor \frac{b^{2}-l_{1}^{2}+l_{1}^{2}+c}{b-l_{1}%
	}\right\rfloor =2b+\left\lfloor \frac{l_{1}^{2}+c}{b-l_{1}}\right\rfloor
	\text{;}$$
	$$f(y_{2})=y_{2}+\left\lfloor \frac{m}{y_{2}}\right\rfloor =b-l_{2}%
	+\left\lfloor \frac{m}{b-l_{2}}\right\rfloor =$$
	$$=b-l_{2}+\left\lfloor \frac{b^{2}-l_{2}^{2}+l_{2}^{2}+c}{b-l_{2}%
	}\right\rfloor =2b+\left\lfloor \frac{l_{2}^{2}+c}{b-l_{2}}\right\rfloor .
	$$
	Now the claim follows from $l_{1}\geq l_{2}$.
	
	b) $x_{m}=b+1$. Let $y_{1}+l_{1}=b+1$, $y_{2}+l_{2}=b+1$, $l_{1},l_{2}\in\mathbb{R}$. It holds $l_{1}\geq l_{2}$. We have%
	$$
	f(y_{1})=y_{1}+\left\lfloor \frac{m}{y_{1}}\right\rfloor =b+1-l_{1}
	+\left\lfloor \frac{m}{b+1-l_{1}}\right\rfloor =$$
	$$=b+1-l_{1}+\left\lfloor \frac{(b+1)^{2}-l_{1}^{2}+l_{1}^{2}+c-2b-1}=
	{b+1-l_{1}}\right\rfloor =$$
	$$=2b+2+\left\lfloor \frac{l_{1}^{2}+c-2b-1}{b+1-l_{1}%
	}\right\rfloor \text{;}$$
	$$f(y_{2})=y_{2}+\left\lfloor \frac{m}{y_{2}}\right\rfloor =b+1-l_{2}%
	+\left\lfloor \frac{m}{b+1-l_{2}}\right\rfloor =$$
	$$=b+1-l_{2}+\left\lfloor \frac{(b+1)^{2}-l_{2}^{2}+l_{2}^{2}+c-2b-1}%
	{b+1-l_{2}}\right\rfloor =$$
	$$=2b+2+\left\lfloor \frac{l_{2}^{2}+c-2b-1}{b+1-l_{2}%
	}\right\rfloor \text{,}%
	$$
	and the claim again follows from $l_{1}\geq l_{2}$.
	
	c) $x_{m}=b+2$. It is easy to see that this case follows analogously as the
	case b).
	
	We have shown $f(y_{1})\geq f(y_{2})$ for $y_{1},y_{2}\in\mathbb{R}$ such that $y_{1}\leq y_{2}\leq x_{m}$, so it holds that $f|_{\left\langle
		-\infty,x_{m}\right\rangle }$ is monotonically decreasing.
\end{proof}

%--------------------------------------------------------------------------------------

\begin{theorem}\label{glavni}
    For $m\in\mathbb{N}$ it holds
	$$
	K(1,m,n)=\left\{
	\begin{array}
	[c]{cc}\scriptstyle
	+\infty, & n\leq2+m+\sqrt{4m+1}\text{;}\\
	\left\lfloor \scriptstyle\frac{1}{2}(n-m-\sqrt{n^{2}+m^{2}-4n-2mn+4})\right\rfloor \scriptstyle+1, &
	n>2+m+\sqrt{4m+1}\text{.}
	\end{array}
	\right.\
	$$
\end{theorem}

\begin{proof}
    Let us define the set $S$ of ordered pairs $(m,n)$ such that there is at least
	one $k$ such that $(1,m,n,k)$ can be realized. We will prove the theorem in
	several steps. First, we will show that for a given $\left(  m,n\right)  \in
	S$, the minimal $k$ for which $(1,m,n,k)$ can be realized, i.e. $k=K(1,m,n)$,
	satisfies the inequality
	$$
	n\geq k+m+\left\lfloor \frac{m}{k-1}\right\rfloor +2\text{.}%
	$$
	Then, by solving that inequality we will obtain the minimal $n$ for which
	$K(1,m,n)\neq\infty$, for given $m\in\mathbb{N}$. Last, we will find the minimal     number of colors needed for given $m$ and $n$.
	
	First we prove two auxiliary claims.
	
	\textbf{CLAIM 1)} Let $(G,\kappa)$ be a pair that realizes $(1,m,n,k)$ for some
	$m,n,k\in\mathbb{N}$ and let $\Delta(G)$ be the maximal degree in $G$. It holds
	$$
	n\geq\Delta(G)+m+\left\lfloor \frac{m}{k-1}\right\rfloor +3\text{.}
	$$

	Let us prove this. Suppose to the contrary, that
	
	$$n\leq\Delta(G)+m+\left\lfloor \frac{m}{k-1}\right\rfloor +2\text{.}$$
 
    Because the $1$-HR condition holds, $\kappa (v)\neq\{1,...,k\}$, for all $v\in V(G)$. So each vertex is missing at least one color, and without the loss of generality we may assume that every vertex is missing exactly one color, Lemma \ref{prosirenje}.
	
	For $u\in V$ such that $d(u)=\Delta(G)$, $G^{\prime}=G-N_{0}(u)$ has
	$p-\Delta(G)-1$ vertices. Let us define the sets $V_{1},...,V_{k}\subseteq
	V(G^{\prime})$ so that a vertex $v$ is in the set $V_{i}$ if and only if $i\notin \kappa(v)$. Every vertex is missing exactly one color so all the vertices in
	$G^{\prime}$ belong in precisely one of $V_{1},...,V_{k}$. Let
	$V_{l}$ be the set for which $\left\vert V_{l}\right\vert \geq\left\vert V_{1}\right\vert
	,...,\left\vert V_{k}\right\vert $ holds. Note that $\left\vert V_{l}\right\vert
	\geq\left\lceil \frac{n-\Delta(G)-1}{k}\right\rceil $. For any set $M$, $|M|=m$, such that $M\supseteq V\left(  G^{\prime}\right)  \backslash
	V_{l}$, all the vertices in $G'-M$ are missing the color $l$ so the $k$-multicoloring of $G$ is not highly $(1,m)$-resistant. This will be possible when
	$$
	n-\Delta(G)-1-\left\lceil \frac{n-\Delta(G)-1}{k}\right\rceil \leq m\text{.}
	$$
	Let us prove that this holds. By the assumption we have
	$$
	n-\Delta(G)-1\leq m+\left\lfloor \frac{m}{k-1}\right\rfloor +1\text{,}
	$$
	and if we write $m=b(k-1)+c$, $0\leq c<k-1$, we have
	$$n-\Delta(G)-1-\left\lceil \frac{n-\Delta(G)-1}{k}\right\rceil\leq$$
	$$\leq
	m+\left\lfloor \frac{m}{k-1}\right\rfloor +1-\left\lceil \frac{m+\left\lfloor
		\frac{m}{k-1}\right\rfloor +1}{k}\right\rceil\leq$$
	$$\leq m+b+1-\left\lceil \frac{b(k-1)+c+b+1}{k}\right\rceil=$$
	$$=m+b+1-b-\left\lceil \frac{c+1}{k}\right\rceil \leq m\text{.}$$
	which proves the claim. We conclude that
	
	$$ p\geq\Delta(G)+m+\left\lfloor \frac{m}{k-1}\right\rfloor +3\text{.}$$

 \textbf{CLAIM 2)} Let $m,n\in\mathbb{N}$ and let $(G,\kappa)$ be a pair which realizes $(1,m,n,k)$, for some $k\in\mathbb{N}$. Let $\Delta(G)$ be the highest degree in $G$. It holds
	$$
	n\geq\Delta(G)+m+\left\lfloor \frac{m}{\Delta(G)}\right\rfloor +3\text{.}
	$$
	Let us prove this by contradiction. Let $u$ be the vertex with the highest degree in $G$. For
	$A=\{u\}$, $G-N_{0}(A)$ has $n-\Delta(G)-1$ vertices. We will show
	that if
	$$
	n\leq\Delta(G)+m+\left\lfloor \frac{m}{\Delta(G)}\right\rfloor +2
	$$
	holds, than there will exist a set $M\subseteq V(G)$, $|M|=m$, such that $G-(N_{0}(A)\cup M)$ is an empty graph. Let us call a vertex pruned if all
	its neighbors are in the set $M$. Since $\Delta(G)$ is the highest degree, there exists a set $M$ such that there are at
	least $\left\lfloor \dfrac{m}{\Delta(G)}\right\rfloor $ pruned vertices in $G-(N_{0}(A)\cup M)$. If
	$G-N_{0}(A)$ has $\left\lfloor \dfrac{m}{\Delta(G)}\right\rfloor
	+m+1$ vertices, such $M$ leads to at most one vertex in $G-(N_{0}(A)\cup M)$, that is neither in $M$ nor a pruned vertex. However than it is an isolated vertex in
	$G-(N_{0}(A)\cup M)$. So we conclude that if
	$$
	n\leq\Delta(G)+m+\left\lfloor \frac{m}{\Delta(G)}\right\rfloor +2
	$$
	holds, no graph $G$ with $n$ vertices can realize $(1,m,p,k)$, for any $k\in\mathbb{N}$ because all the vertices in $G-(N_{0}(A)\cup M)$ will be isolated, and the $1$-HR condition must hold. Therefore
	$$
	n\geq\Delta(G)+m+\left\lfloor \frac{m}{\Delta(G)}\right\rfloor +3\text{.}
	$$

	This concludes the proof of Claim 2.

	Let $m\in\mathbb{N}$ and let
	$$
	n_{m}=\min_{\Delta_{0}\in\mathbb{N}}\{\Delta_{0}+m+\left\lfloor \frac{m}{\Delta_{0}}\right\rfloor +3\}\text{,}
	$$
	and let $\Delta_{m}$ be a value of $\Delta_{0}$ for which this minimum is
	obtained. (If there is more than one such value choose the smallest one).
	
	Furthermore, let us prove that $(1,m,n_{m},\Delta_{m}+1)$ can be realized. The
	realization is a pair $(G_{m},\kappa)$ given as follows
	
	$$
	G_{m}=\underset{\left\lfloor \dfrac{n_{m}}{\Delta_{m}+1}\right\rfloor \text{
			times}}{\underbrace{K_{\Delta_{m}+1}\cup K_{\Delta_{m}+1}\cup...\cup
			K_{\Delta_{m}+1}}}\cup K_{l},$$
	$$\text{ where }l=n_{m}-\left\lfloor \dfrac{n_{m}%
	}{\Delta_{m}+1}\right\rfloor \cdot\left(  \Delta_{m}+1\right)  \text{,}
	$$
	with function $\kappa$ that gives each vertex all but one color in such a way that
	every two vertices in the same component miss different colors. The realization
	of some examples is illustrated in Figure \ref{fig:1_11_16}.
	
	\begin{figure}[h]
		\centering\includegraphics[scale=0.6]{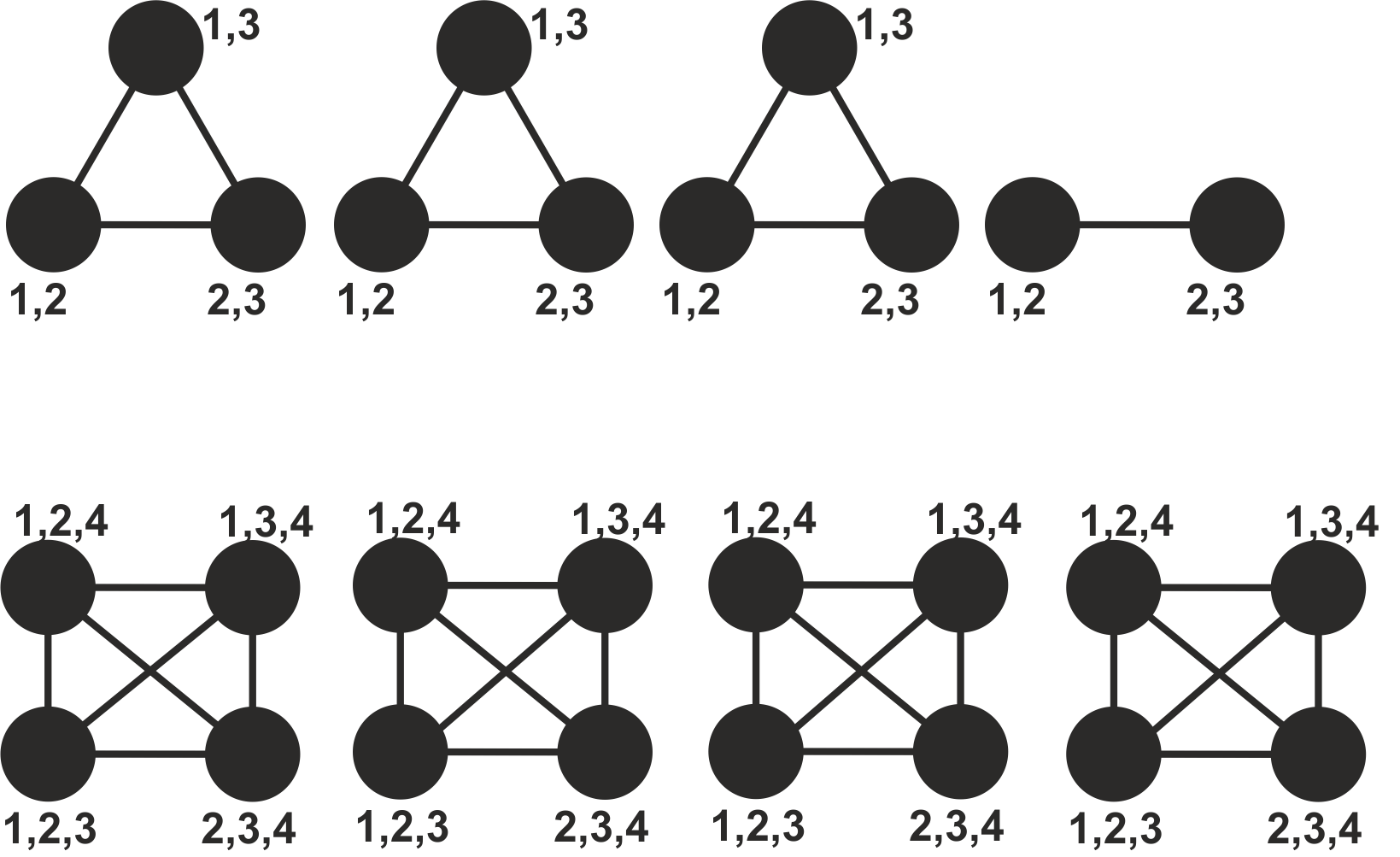}
		\caption{Realizations for $(1,4,11,3)$ and $(1,8,16,4)$}
		\label{fig:1_11_16} \ \
	\end{figure}

 We will now prove that the smallest $k$, for which $(1,m,n,k)$ can be realized
	for given $\left(  m,n\right)  \in S$, satisfies the expression
	
	$$n\geq k+m+\left\lfloor \frac{m}{k-1}\right\rfloor +2\text{.}$$
	Let $\left(  m,n\right)  \in S$ and let $k_{0}$ be the smallest $k$ for which
	$(1,m,n,k)$ can be realized. We consider 3 cases.
	
	1) $k_{0}>\Delta_{m}+1$.
	
	From Claim 2 we have
	
	$$n\geq\Delta(G)+m+\left\lfloor \frac{m}{\Delta(G)}\right\rfloor +3\geq$$
	$$\geq
	\Delta_{m}+m+\left\lfloor \frac{m}{\Delta_{m}}\right\rfloor +3\text{.}$$

	We know\ that $(1,m,n_{m},\Delta_{m}+1)$ can be realized so because $n\geq
	n_{m}$, from Proposition \ref{uvodna} it follows that $(1,m,n,\Delta_{m}+1)$ can be
	realized. But now obviously the chosen $k_{0}$ is not the smallest $k$ for
	which $(1,m,n,k)$ can be realized, so this case cannot happen.
	
	2) $k_{0}=\Delta_{m}+1$. It now holds
	$$
	n\geq\Delta_{m}+m+\left\lfloor \frac{m}{\Delta_{m}}\right\rfloor
	+3=k_{0}+m+\left\lfloor \frac{m}{k_{0}-1}\right\rfloor +2\text{,}
	$$
	so the claim stands.
	
	3) $k_{0}<\Delta_{m}+1$. Let $G$ be a graph which realizes $(1,m,n,k_{0})$,
	and let $\Delta(G)$ be the highest degree in $G$. We observe 2 subcases.
	
	3.1.) $k_{0}\leq\Delta(G)+1$. Now from Claim 1 it holds
	$$
	n\geq\Delta(G)+m+\left\lfloor \frac{m}{k_{0}-1}\right\rfloor +3\geq
	k_{0}+m+\left\lfloor \frac{m}{k_{0}-1}\right\rfloor +2\text{,}
	$$
	so the claim stands.
	
	3.2.) $k_{0}>\Delta(G)+1$. Now we have
	$$
	\Delta(G)+1<k_{0}<\Delta_{m}+1\text{,}
	$$
	so from Claim 2 it follows
	$$
	n\geq\Delta(G)+m+\left\lfloor \frac{m}{\Delta(G)}\right\rfloor +3\text{,}
	$$
	and from Lemma \ref{pad} it follows that
	$$
	n\geq\Delta(G)+m+\left\lfloor \frac{m}{\Delta(G)}\right\rfloor +3\geq
	k_{0}+m+\left\lfloor \frac{m}{k_{0}-1}\right\rfloor +2\text{,}$$
	because $\Delta(G)<k_{0}\leq\Delta_{m}$, and the function $f(\Delta
	(G))=\Delta(G)+1+m+\left\lfloor \frac{m}{\Delta(G)}\right\rfloor +2$ is
	monotonically decreasing on $f|_{\left\langle -\infty,\Delta_{m}\right\rangle
	}$.
	
	We have proven that the minimal number of colors, $k_{0}\geq2$, such that
	$(1,m,n,k_{0})$ can be realized satisfies
	
	\begin{equation}\label{eq:1}
	\centering
	n\geq k_{0}+m+\left\lfloor \frac{m}{k_{0}-1}\right\rfloor +2\text{, }
	\end{equation}
	
	for given $\left(  m,n\right)  \in S$. Moreover, graph
	$$
	G_{0}=\underset{\left\lfloor \dfrac{n}{k_{0}}\right\rfloor \text{ times}
	}{\underbrace{K_{k_{0}}\cup K_{k_{0}}\cup...\cup K_{k_{0}}}}\cup K_{l},\text{
		where }l=n-\left\lfloor \dfrac{n}{k_{0}}\right\rfloor \cdot k_{0}\text{,}
	$$
	with function $\kappa$ that gives each vertex all but one color, in such a way that
	every two vertices in the same component miss different colors, realizes
	$(1,m,n,k_{0})$. Hence indeed
	$$
	K(1,m,n)=k_{0}\text{.}
	$$

 Further, for every $m,n\in\mathbb{N}$, if there is $k_{0}$ that satisfies (1), then $(m,n)\in S$. Moreover, (1)
	implies that for each $m\in\mathbb{N}$ there exists an $n\in\mathbb{N}$ such that $(m,n)\in S$.
	
	Now let us find $K(1,m,n)$, and the minimal $n$ for which $K(1,m,n)\neq
	\infty$.
	
	Let $m,n\in\mathbb{N}$. Note that right hand side of inequality (1) is an integer, hence (1) is
	equivalent to
	$$
	n+1>k_{0}+m+\frac{m}{k_{0}-1}+2\text{.}
	$$
	This inequality is equivalent to
	$$
	k_{0}^{2}+k_{0}(m-n)+n-1<0\text{,}
	$$
	and by solving it we get the interval
	$$
	k_{0}\in\left\langle \frac{1}{2}(n-m-\sqrt{D}),\frac{1}{2}(n-m+\sqrt
	{D})\right\rangle \text{,}
	$$
	where $D=n^{2}+m^{2}-4n-2mn+4$. Let us denote that interval by $I_{k}$. So
	$K(1,m,n)\neq\infty$ if and only if there exists an integer $k_{0}>1$ in $I_{k}
	$. Let us find the conditions for $n$, for given $m\in\mathbb{N}$.
	
	Obviously, if $D<0$ than $K(1,m,n)=\infty$ because the value for $k_{0}
	^{2}+k_{0}(m-n)+n-1$ is never negative. For $D=0$ we obtain $I_{k}
	=\emptyset$, and it remains to observe the case $D>0$. The values for $n$, for
	which $D>0$, are in the set
	$$
	\left(  \left\langle -\infty,2+m-2\sqrt{m}\right\rangle \cup\left\langle
	2+m+2\sqrt{m},+\infty\right\rangle \right)  \cap\mathbb{N}\text{.}
	$$
	Let us define functions $n_{1},n_{2}:\mathbb{N}\rightarrow\mathbb{R}$ with
	
	$$n_{1}(m)=2+m-2\sqrt{m},$$
	$$n_{2}(m)=2+m+2\sqrt{m}.$$
	
	$n<n_{1}(m)$ would imply that there exists $k\in\mathbb{N}$ such that $(1,m,n,k)$ can be realized. Then, because of Proposition \ref{uvodna}, such $k$ could
	also be found for all $n\in\mathbb{N}$, $n\geq n_{1}(m)$ and we have already proven that is impossible for
	$n\in\left[  n_{1}(m),n_{2}(m)\right]  \cap\mathbb{N}$ (note that this set is non-empty). So we conclude that $n>n_{2}(m)$ must
	hold. Hence,
	$$
	k_{0}=\left\lfloor \frac{1}{2}(n-m-\sqrt{D})\right\rfloor +1,
	$$
	and
	$$
	\left\lfloor \frac{1}{2}(n-m-\sqrt{D})\right\rfloor +1<\frac{1}{2}%
	(n-m+\sqrt{D}),\text{ and }n>n_{2}(m)\text{.}%
	$$
	Let us distinguish two cases.
	
	CASE 1) $\frac{1}{2}(n-m-\sqrt{D})+1<\frac{1}{2}(n-m+\sqrt{D}).$ 
	\\
	This is
	equivalent to $D>1$ which is further equivalent to
	$$
	n\in\left(  \left\langle -\infty,2+m-\sqrt{4m-1}\right\rangle \cup\left\langle
	2+m+\sqrt{4m+1},+\infty\right\rangle \right)  \cap\mathbb{N}\text{.}%
	$$
	Let us define functions $n_{3},n_{4}:\mathbb{N}\rightarrow\mathbb{R}$ with%
	
	$$n_{3}(m)=2+m-\sqrt{4m-1}\text{,}$$
	$$n_{4}(m)=2+m+\sqrt{4m+1}\text{.}$$
	
	Since $n_{3}(m)<n_{1}(m)<n_{2}(m)<n_{4}(m)$, for all $m\in\mathbb{N}$, this reduces to $n>n_{4}(m)$.
	
	CASE 2) $\frac{1}{2}(n-m-\sqrt{D})+1\geq\frac{1}{2}(n-m+\sqrt{D})$. \\
	Then
	$n_{2}(m)<n\leq n_{4}(m)$.
	
	We will now prove that $I_{k}\cap\mathbb{N}=\emptyset$ holds. Let us distinguish two subcases.
	
	2.1) $n<n_{4}(m)$.
	
	We need to prove that the smallest integer larger than $n_{2}\left(  m\right)
	$ is larger than $n_{4}\left(  m\right)  ,$ i.e. that
	
	$$\left\lfloor n_{2}\left(  m\right)  \right\rfloor +1\geq n_{4}\left(
	m\right);$$
	$$\left\lfloor 2+m+2\sqrt{m}\right\rfloor +1\geq2+m+\sqrt{4m+1};$$
	$$1+\left\lfloor 2\sqrt{m}\right\rfloor\geq\sqrt{4m+1}.$$

	Let us denote $\left\lfloor 2\sqrt{m}\right\rfloor =b$. We may write
	$4m=b^{2}+c$, $0\leq c\leq2b$. Now we have to prove
	$$
	1+b\geq\sqrt{b^{2}+c+1}.
	$$
	However, further calculation gives us that this holds when $2b\geq c$, so the
	claim is proven.
	
	2.2) For $n=n_{4}(m)$ it holds $I_{k}\cap\mathbb{N}=\emptyset$. Let us observe the case $n=n_{4}(m)=2+m+\sqrt{4m+1}$. Since
	$n\in\mathbb{N}$, this is possible only if $\sqrt{4m+1}\in\mathbb{N}$. For this value of $n$, $I_{k}$ is equal to%
	$$
	I_{k}=\left\langle \frac{1}{2}+\frac{1}{2}\sqrt{4m+1},\frac{3}{2}+\frac{1}%
	{2}\sqrt{4m+1}\right\rangle \text{.}%
	$$
	It is easy to see that now the length of $I_{k}$ is precisely $1$, and that the
	bounds of the interval are natural numbers. But this means $I_{k}\cap\mathbb{N}=\emptyset$.
	
	We have proven that no $n\in\left\langle n_{2}(m),n_{4}(m)\right]  \cap\mathbb{N}$ can produce $k\in I_{k}\cap\mathbb{N}$ so we conclude that $K(1,m,n)=\infty$ if and only if $n\leq n_{4}(m)$ and
	for $n>n_{4}(m)$ we have
	$$
	K(1,m,n)=\left\lfloor \frac{1}{2}(n-m-\sqrt{n^{2}+m^{2}-4n-2mn+4}
	)\right\rfloor +1\text{,}
	$$
	which concludes the proof of the theorem. 
\end{proof}

\begin{remark}
    Theorem \ref{glavni} solves the problem of finding the minimal number of nodes and the minimal number of pieces a secret must be divided to, in order for the network to be resilient to the attack of $1$ agent that steals the pieces given to him and betrays its neighbors, and $m\in\mathbb{N}$ missing or malfunctioning nodes. It also gives one possible network configuration and the pieces distribution. For example, if there are $8$ malfunctioning nodes, as well as $1$ agent, the network must have the minimum of $16$ nodes and the secret has to be divided into minimum of $4$ parts in order to preserve it. One possible distribution is given in Figure \ref{fig:1_11_16}.
\end{remark}

\begin{remark}
    The case of $3$ attackers and $1$ malfunctioning node is analyzed in \cite{nas2}.
\end{remark}

%-------------------------------------------------------------

\end{document}